\documentclass[preprint,12pt]{elsarticle}
\newtheorem{thm}{Theorem}[section]
\newtheorem{corollary}{Corollary}[section]
\newtheorem{lemma}{Lemma}[section]
\newtheorem{remark}{Remark}[section]

\usepackage{amssymb}

\begin{document}

\begin{frontmatter}

\title{Hoeffding's inequality for supermartingales}
\author{Xiequan Fan$^*$,\ \ \ \ \ \ Ion Grama\ \ \  \ \ and\ \ \ \ \ \ Quansheng Liu}
 \cortext[cor1]{\noindent Corresponding author. \\
\mbox{\ \ \ \ }\textit{E-mail}: fanxiequan@hotmail.com (X. Fan), \ \ \ \ \  ion.grama@univ-ubs.fr (I. Grama),\\
\mbox{\ \ \ \ \ \ \ \ \ \ \  \ \ \ \ }quansheng.liu@univ-ubs.fr (Q. Liu). }
\address{Universit\'{e} de Bretagne-Sud, LMBA, UMR CNRS 6205, \\
 Campus de Tohannic, 56017 Vannes, France}

\begin{abstract}

We give an extension of  Hoeffding's inequality to the case of
supermartingales with differences bounded from above. Our inequality strengthens or extends
the inequalities of Freedman, Bernstein, Prohorov, Bennett and  Nagaev.
\end{abstract}

\begin{keyword} Concentration inequalities;
Hoeffding's inequality; Freedman's inequality; Bennett's inequality; martingales; supermartingales
\vspace{0.3cm}
\MSC Primary 60G42 \sep 60G40 \sep 60F10; second 60E15 \sep 60G50

\end{keyword}

\end{frontmatter}


\section{Introduction}
Let $\xi _{1},...,\xi _{n}$ be a sequence of centered ($\mathbb{E}\xi_{i}=0$) random
variables such that $\sigma _{i}^{2}=\mathbb{E}\xi _{i}^{2}<\infty $ and let $
X_{n}=\sum_{i=1}^{n}\xi _{i}.$ Since the seminal papers of Cram\'{e}r (1938, \cite{Cramer38}) and Bernstein
(1946, \cite{B46}), the estimation of the
tail probabilities $\mathbb{P}\left( X_{n}>x \right)$ for positive $x$ has attracted
much attention. We would like to mention here the celebrated Bennett inequality
(1962, cf. (8b) of \cite{B62},  see also Hoeffding \cite{Ho63}) which states that, for independent
and centered random variables $\xi_{i}$ satisfying $\xi _{i}\leq 1$ and for any $t>0,$%
\begin{eqnarray}
\mathbb{P}(X_{n}\geq nt )&\leq& \left( \frac{ \sigma^2 }{t+ \sigma^2}\right) ^{n(t+ \sigma^2)}e^{nt} \label{Bennett ineq} \\
 &=& \exp\left\{ -n t \left[ \left(1+\frac{\sigma^2}{t}\right)\log\left(1+\frac{t}{\sigma^2}\right)-1 \right]\right\} ,
\end{eqnarray}%
where $\sigma ^{2}=\frac{1}{n}\sum_{i=1}^{n}\sigma _{i}^{2}$. Further, inequalities for the probabilities
$\mathbb{P}\left( X_{n}>x\right) $ have been obtained by Prohorov (1959, \cite{Pr59}%
), Hoeffding (1963, \cite{Ho63}), Azuma (1967, \cite{A67}), Steiger (1967, \cite{St67}; 1969, \cite{St69}), Freedman (1975,
\cite{FR75}), Nagaev (1979, \cite{N79}), Haeusler (1984, \cite{H84}), McDiarmid (1989, \cite{M}),
Pinelis (1994, \cite{P94a}), Talagrand (1995, \cite{Ta95}), De La Pe\~{n}a
(1999, \cite{D99}), Lesigne and Voln\'{y} (2001, \cite{LV01}), Nagaev (2003,
\cite{N03}), Bentkus (2004, \cite{Be04}), Pinelis (2006, \cite{P06})
and Bercu and Touati (2008, \cite{BT08}) among others.

Most of these results were obtained by an approach based
on the use of the exponential Markov's inequality. The challenge for this
method is to find a sharp upper bound of the moment generating function
 $\phi _{i}\left( \lambda \right) =\mathbb{E}(e^{\lambda \xi _{i}}).$
Hoeffding \cite{Ho63}, Azuma \cite{A67} and McDiarmid \cite{M} used the
elementary estimation $\phi _{i}\left( \lambda \right) \leq e^{\lambda
^{2}/2},\lambda \geq 0,$ which holds if $|\xi _{i}|\leq 1.$ Better results can be obtained by the following
improvement $\phi _{i}\left( \lambda \right) \leq (e^{\lambda }-1-\lambda
)\sigma _{i}^{2},\ \lambda \geq 0,$ which holds for $\xi _{i}\leq 1$
(see for example Freedman \cite{FR75}). Bennett \cite{B62} and Hoeffding \cite{Ho63} used a more precise estimation
\begin{equation}
\phi _{i}\left( \lambda \right) \leq \frac{1}{1+\sigma _{i}^{2}}\exp \left\{
-\lambda \sigma _{i}^{2}\right\} +\frac{\sigma _{i}^{2}}{1+\sigma _{i}^{2}}%
\exp \{\lambda \},\ \ \ \ \lambda \geq 0,  \label{be}
\end{equation}%
for any $\xi_{i}$ satisfying $\xi _{i}\leq 1.$  Bennett's estimation (\ref{be}) is sharp with the equality  attained when $\mathbb{P}(\xi _{i}=1)=\frac{%
\sigma _{i}^{2}}{1+\sigma _{i}^{2}}$ and$\ \mathbb{P}(\xi _{i}=-\sigma _{i}^{2})=%
\frac{1}{1+\sigma _{i}^{2}}.$

 Using (\ref{be}), Hoeffding improved Bennett's inequality
(\ref{Bennett ineq}) and obtained the following inequality:
for independent and centered random variables $(\xi _i)_{i=1,...,n}$ satisfying $\xi_i\leq 1$ and for any $0< t < 1$,
\begin{eqnarray}
\mathbb{P}\left(X_n\geq nt \right) &\leq& \left\{ \left( 1+ \frac{t}{ \sigma^2}\right) ^{-\frac{t+\sigma^2}{1+ \sigma^2} }\left(1-t\right) ^{-\frac{1-t}{1+ \sigma^2} }\right\} ^n , \label{Ho}
\end{eqnarray}
where $\sigma^{2}=\frac{1}{n}\sum_{i=1}^{n}\sigma _{i}^{2}$ (cf. (2.8) of \cite{Ho63}).

It turns out that, under the
stated conditions, Hoeffding's inequality (\ref{Ho}) is very tight and improving (\ref{Ho}) is a rather difficult task.
Significant advances in improving Hoeffding's and Bennett's inequalities have been obtained by several
authors. For instance Eaton \cite{E74}, Pinelis \cite{P94} and Talagrand \cite{Ta95}
have added to (\ref{Ho}) a missing factor of the order $\frac{\sigma }{\sqrt{n}\,\, t}$.
Improvements of the Bennett's inequality (\ref{Bennett ineq}) can be found in Pinelis
\cite{P06}, where some larger classes of functions are
considered instead of  the class of exponential functions usually used in
Markov's inequality. When $\xi _{i}$ are martingale differences,
\mbox{Bentkus} \cite{Be04} showed that if the conditional variances of $\xi _{i}$ are
bounded, then $\mathbb{P}(X_n\geq x)\leq c\ \mathbb{P}(\sum_{i=1}^{n}\eta
_{i}\geq x),$ where $\eta _{i}$ are independent and identically
distributed  Rademacher random variables, $c=e^{2}/2=3.694...$~and $x$ is a
real such that $\mathbb{P}(\sum_{i=1}^{n}\eta _{i}\geq x)$ has a jump at $x$ (see
also \cite{Be06} for related results). However, to the best
of our knowledge, there is no martingale or supermartingale version which
reduces exactly to the Hoeffding inequality (\ref{Ho}) in the independent
case.

The scope of the paper is to extend the Hoeffding inequality (\ref{Ho}) to
the case of martingales and supermartingales. Our inequality will recover  (\ref{Ho})
in the independent case, and in the case of (super)martingales
will apply under a very weak  constraint  on the sum of conditional variances.

The main results of the paper are the following inequalities (see Theorem %
\ref{th1} and Remark \ref{remark1}). Assume that $(\xi _i,\mathcal{F}_i)_{i=1,...,n}$ are supermartingale
differences satisfying $\xi _{i}\leq 1$. Denote by $\langle X\rangle
_{k}=\sum_{i=1}^{k}\mathbb{E}(\xi _{i}^{2}|\mathcal{F}_{i-1})$ for $k=1,...,n$. Then,
for any $x\geq 0$ and $v>0$,
\begin{eqnarray}
&& \mathbb{P}\left( X_k \geq x\ \mbox{and}\ \langle X\rangle_{k}\leq v^2\ \mbox{for some}\ k \in [1, n] \right) \nonumber \\
&\leq& \left\{\left( \frac{v^2}{
x+v^2}\right)^{ x+v^2 }\left( \frac{n}{n-x}\right)^{ n-x }
\right\}^{\frac{n}{n+v^2} }\mathbf{1}_{\{x \leq n\}} \label{fglf}\\
&\leq& \exp\left\{-\frac{x^2}{2(v^2+ \frac{1}{3}x )}\right\}.\label{fbem}\ \ \ \ \ \ \ \ \ \ \ \ \ \ \ \ \ \ \ \ \ \ \ \ \ \ \ \ \ \  \ \ \ \
\end{eqnarray}
In the independent case, inequality (\ref{fglf}) with $x=nt$ and $v^2=n\sigma^2$ reduces to inequality (\ref{Ho}).
We will see that the inequalities (\ref{fglf}) and (\ref{fbem}) strengthen or extend many well-known inequalities obtained by Freedman, De La Pe\~{n}a,
Bernstein, Prohorov, Bennett, Hoeffding, Azuma,  Nagaev and Haeusler.
In particular, if the  martingale differences $(\xi _i,\mathcal{F}_i)_{i=1,...,n}$ satisfy $-b \leq \xi_{i} \leq 1$ for some constant $  b >0 $,
then we get (see Corollary \ref{co2}), for all $x\geq0$,
\begin{eqnarray}
\mathbb{P}\left( \max_{ 1 \leq k \leq n} X_{k} \geq x \right)&\leq&   \exp\left\{ - \frac{ 2x^2}{ U_n (x,b) }\right\} ,  \label{f6}
\end{eqnarray}
where
$$U_n(x,b)=\min\left\{n(1+b)^2,\  4\left(nb+\frac {1}{3} x \right)\right\}.$$
Notice that inequality (\ref{f6}) is sharper than the usual Azuma-Hoeffding  inequality when $0< x < \frac{3}{4}n(1-b)^2$.

Our approach is based on the conjugate distribution technique due to Cram\'{e}r, and is different from the method used in
Hoeffding's original paper \cite{Ho63}. The technique has been developed in Grama and Haeusler \cite{GH00}
to obtain expansions of large deviation for martingales. We refine this technique
to get precise upper bounds for tail
probabilities, providing a simple and unified approach for improving several well-known inequalities. We also
make clear some relations among these inequalities.

Our main results will be presented in Section \ref{sec2} and proved in Sections \ref{sec3} and \ref{sec4}.

\section{Main Results} \label{sec2}
Assume that we are given a sequence of real supermartingale differences $%
(\xi _i,\mathcal{F}_i)_{i=0,...,n}, $ defined on some probability space $%
(\Omega ,\mathcal{F},\mathbb{P})$, where $\xi _0=0 $ and $\{\emptyset, \Omega\}=%
\mathcal{F}_0\subseteq ...\subseteq \mathcal{F}_n\subseteq \mathcal{F}$ are
increasing $\sigma$-fields. So, by definition, we have $\mathbb{E}(\xi_{i}|\mathcal{F}_{i-1})\leq 0, \  i=1,...,n$. Set
\begin{equation}  \label{matingal}
X_{0}=0,\ \ \ \ \ X_k=\sum_{i=1}^k\xi _i,\quad k=1,...,n.
\end{equation}
Let $\left\langle X\right\rangle $ be the quadratic characteristic of the
supermartingale $X=(X_k,\mathcal{F}_k)$:
\begin{equation}  \label{quad}
\left\langle X\right\rangle _0=0,\ \ \ \ \ \left\langle X\right\rangle
_k=\sum_{i=1}^k\mathbb{E}(\xi _i^2|\mathcal{F} _{i-1}),\quad k=1,...,n.
\end{equation}

\begin{thm}
\label{th1} Assume that $(\xi _i,\mathcal{F}_i)_{i=1,...,n}$ are supermartingale differences satisfying $\xi_{i} \leq 1$. Then,
for any $ x \geq 0$ and $v > 0$,
\begin{eqnarray}
  \mathbb{P}\left( X_k \geq x\ \mbox{and}\ \langle X\rangle_{k}\leq v^2\ \mbox{for some}\ k \in [1, n] \right) \leq H_{n}(x, v),\label{fgl1}
\end{eqnarray}
where
\[
H_{n}(x,v)=\left\{\left( \frac{v^2}{
x+v^2}\right)^{ x+v^2 }\left( \frac{n}{n-x}\right)^{ n-x }
\right\}^{\frac{n}{n+v^2} }\mathbf{1}_{\{x \leq n\}}
\]
with the convention that $(+\infty)^0=1$ (which applies when $x=n$).
\end{thm}

Because of the obvious inequalities
\begin{eqnarray}
&&\mathbb{P}\left(  X_n \geq x, \langle X\rangle_{n}\leq v^2 \right) \label{Proba1}\\
&\leq& \mathbb{P}\left( \max_{ 1 \leq k \leq n} X_k \geq x, \langle X\rangle_{n}\leq v^2 \right) \label{Proba2}\\
&\leq&\mathbb{P}\left( X_k \geq x\ \mbox{and}\ \langle X\rangle_{k}\leq v^2\ \mbox{for some}\ k \in [1,n] \right),\nonumber
\end{eqnarray}
the function $H_{n}(x,v)$ is also an upper bound of the tail probabilities  (\ref{Proba1})  and (\ref{Proba2}).
Therefore Theorem \ref{th1} extends Hoeffding's inequality (\ref{Ho}) to the case of supermartingales with differences $\xi_i$ satisfying $\xi_i\leq1$.

The following remark establishes some relations among the well-known bounds of Hoeffding, Freedman, Bennett, Bernstein and De La Pe\~{n}a.
\begin{remark}\label{remark1}
  For any $x\geq0$ and $v>0$, it holds
\begin{eqnarray}
H_{n}(x, v)&\leq& F(x,v)=:\left(\frac{v^2}{x+v^2} \right)^{x+v^2}e^x  \label{freedma1} \\
&\leq& B_1(x,v)=:\exp\left\{-\frac{x^2}{v^2\left(1+\sqrt{1+\frac{2\, x}{3\, v^2}  } \right)+\frac 1 3x }\right\} \label{Bennett} \\
&\leq& B_2(x,v)=: \exp\left\{-\frac{x^2}{2(v^2+\frac{1}{3}x )}\right\} . \label{Bernstein}
\end{eqnarray}
Moreover, for any $x, v>0$, $H_n(x,v)$ is increasing in $n$ and
\begin{eqnarray}
 \lim_{n\rightarrow \infty} H_{n}(x, v)  = F(x,v). \label{limn}
\end{eqnarray}
\end{remark}

Since $H_{n}(x,v)\leq F(x,v)$, our inequality (\ref{fgl1}) implies
 Freedman's inequality for supermartingales \cite{FR75}. The bounds $B_1(x,v)$ and $B_2(x,v)$ are
respectively the bounds of Bennett and Bernstein (cf.$\,$ \cite{B62}, (8a) and \cite{B46}). Note that Bennett
and Bernstein obtained their bounds  for independent random variables under the Bernstein condition
\begin{eqnarray}
\mathbb{E} |\xi_{i}|^k  \leq \frac{1}{2}\left(\frac{1}{3}\right)^{k-2}\mathbb{E} \xi_{i}^2,\ \ \  \mbox{for} \ \ \ k\geq3. \label{Bercd}
\end{eqnarray}
We would like to point out that our condition $\xi_i\leq 1$ does not imply  Bernstein condition (\ref{Bercd}).
The bounds $B_1(x,v)$ and $B_2(x,v)$ have also been obtained by De La Pe\~{n}a (\cite{D99}, (1.2)) for martingale differences $\xi_{i}$ satisfying the conditional version of Bernstein's condition  (\ref{Bercd}).
Our result  shows that the inequalities of Bennett (\cite{B62}, (8a)), Bernstein \cite{B46} and De La Pe\~{n}a (\cite{D99}, (1.2)) also hold when the (conditional)
Bernstein condition is replaced by the condition  $\xi_{i}\leq1$.
So Theorem \ref{th1} refines and completes the inequalities of Bennett, Bernstein and
De La Pe\~{n}a for supermartingales with differences bounded from above.

It is interesting to note that from Theorem \ref{th1} and (\ref{Bennett}) it follows that
\begin{eqnarray}\label{beni}
\mathbb{P}\left( X_k\geq \frac x 3+ v\sqrt{2 x} \ \mbox{and}\ \langle X\rangle_k\leq v^2\ \mbox{for some }\ k \in [1,n]\right)\leq e^{-x},
\end{eqnarray}
which is another form  of Bennett's inequalities (for related bounds we refer to Rio \cite{R02} and Bousquet \cite{B02}).

If the (super)martingale differences $(\xi_i,\mathcal{F}_i)_{i=1,..,n}$ are in addition bounded from below, our inequality (\ref{fgl1})
also implies the inequalities (2.1) and (2.6) of Hoeffding \cite{Ho63} as seen from the following corollary.
\begin{corollary} \label{co2}
 Assume that $(\xi _i,\mathcal{F}_i)_{i=1,...,n}$  are  martingale differences satisfying $-b \leq \xi_{i} \leq 1 $ for some constant $b>0$. Then, for any $ x \geq 0$,
\begin{eqnarray}
   \mathbb{P}\left( \max_{ 1 \leq k \leq n} X_k \geq x \right)  \leq   H_{n}\left(x,\sqrt{nb} \right)   \label{Ho11}
\end{eqnarray}
and
\begin{eqnarray}
\mathbb{P}\left( \max_{ 1 \leq k \leq n} X_k \geq x \right)&\leq&   \exp\left\{ - \frac{ 2x^2}{ U_n (x,b) }\right\} ,   \label{Ho12}
\end{eqnarray}
where
$$U_n(x,b)=\min\left\{n(1+b)^2,\  4\left(nb+\frac {1}{3} x \right)\right\}.$$
The inequalities (\ref{Ho11}) and (\ref{Ho12}) remain true for supermartingale differences $(\xi _i,\mathcal{F}_i)_{i=1,...,n}$ satisfying
$-b\leq \xi_i \leq 1$ for some constant $0<b\leq1$.
\end{corollary}

In the martingale case, our inequality (\ref{Ho11}) is a refined version of the
inequality (2.1) of Hoeffding \cite{Ho63} in the sense that  $X_n$ is replaced by $\max_{1\leq k \leq n}X_k$. When $U_{n}(x,b)=n(1+b)^2$, inequality (\ref{Ho12}) is a refined version of the usual Azuma-Hoeffding inequality (cf. \cite{Ho63}, (2.6)); when $0< x < \frac{3}{4}n(1-b)^2$, our inequality (\ref{Ho12}) is sharper than the Azuma-Hoeffding inequality. Related results can be found in Steiger \cite{St67}, \cite{St69}, McDiarmid \cite{M}, Pinelis \cite{P06} and Bentkus \cite{Be04}, \cite{Be06}.

The following result extends an inequality of De La Pe\~{n}a (\cite{D99}, (1.15)).
\begin{corollary} \label{co1}
 Assume that $(\xi _i,\mathcal{F}_i)_{i=1,...,n}$ are supermartingale differences satisfying $ \xi_{i} \leq 1 $. Then, for any $ x\geq 0$ and $v > 0$,
\begin{eqnarray}\label{pie}
&& \mathbb{P}\left( X_k \geq  x\ \mbox{and}\ \langle X\rangle_{k}\leq  v^2 \  \mbox{for some} \ k \in [1,n] \right) \nonumber \\ &\leq& \exp \left\{-\frac{ x}{2}\
   \emph{arc}\sinh \left(\frac{x}{2v^2}\right)\right\} .
\end{eqnarray}
\end{corollary}

De La Pe\~{n}a \cite{D99} obtained the same inequality (\ref{pie}) for martingale differences $(\xi _i,\mathcal{F}_i)_{i=1,...,n}$ under the more
restrictive condition that $|\xi_{i}|\leq c$  for some constant $0< c < \infty$.
In the independent case, the bound in (\ref{pie}) is the Prohorov bound \cite{Pr59}. As was remarked by Hoeffding \cite{Ho63},
the right side of (\ref{fgl1}) is less than the right side of (\ref{pie}). So inequality (\ref{fgl1}) implies inequality (\ref{pie}).

For unbounded supermartingale differences $(\xi _i,\mathcal{F}_i)_{i=1,...,n}$, we have the following inequality.

\begin{corollary}\label{th2.5}
Assume that $(\xi _i,\mathcal{F}_i)_{i=1,...,n}$  are supermartingale differences. Let $y > 0$ and
\begin{equation} \label{bnk}
 V_k^2(y)=\sum_{i=1}^k \mathbb{E}(\xi_{i}^2 \mathbf{1}_{\{\xi_{i}\leq y \}} | \mathcal{F}_{i-1}),\ \ k=1,...,n.
\end{equation}
Then, for any $ x\geq 0,  y> 0$ and $v> 0$,
\begin{eqnarray}
&&\mathbb{P}\left( X_k \geq x\ \mbox{and}\ V_k^2(y)\leq v^2 \ \mbox{for some}\ k \in [1,n] \right) \nonumber \\
&\leq&H_n\left(\frac{x}{y},\frac{v}{y}\right)  +\mathbb{P}\left( \max_{1\leq i \leq n}\xi_{i}> y \right) \label{ffn} .
\end{eqnarray}
\end{corollary}

We notice that inequality (\ref{ffn}) improves an inequality of Fuk (\cite{Fu73}, (3)).
It also extends and improves Nagaev's inequality (\cite{N79}, (1.55))  which was
obtained in the independent case.

Since $\mathbb{P}(V_n^2(y)>  v^2) \leq P(\langle X\rangle_n >
v^2)$ and $H_n(x,v)\leq F(x,v)$, Corollary \ref{th2.5} implies the following inequality due to Courbot \cite{Co99}:
\begin{eqnarray}
 \mathbb{P}\left(\max_{ 1 \leq k \leq n} X_k \geq x \right)
&\leq & F\left(\frac xy, \frac vy \right)   + \sum_{i=1}^{n}\mathbb{P}\left( \xi_{i}> y \right)  +  \mathbb{P}(\langle X\rangle_n >
v^2) .\ \ \ \ \ \
\end{eqnarray}
 A slightly weaker inequality was obtained earlier by Haeusler \cite{H84}:
in Haeusler's inequality $F\left(\frac xy, \frac vy \right)$ is replaced by a larger bound $\exp\left\{\frac{x}{y}\left(1-\log\frac{xy}{v^2}\right)\right\}$. Thus, inequality
(\ref{ffn}) improves Courbot's and Haeusler's inequalities.

To  close this section, we present an  extension of the inequalities of Freedman and Bennett under the condition that
\begin{eqnarray}
\mathbb{E}(\xi_i^2e^{\lambda \xi_{i}}|\mathcal{F}_{i-1}) \leq e^{\lambda}E(\xi_i^2|\mathcal{F}_{i-1}),\ \ \ \ \mbox{for any} \ \ \lambda\geq 0,\label{cdwe}
\end{eqnarray}
which is weaker than the assumption $\xi_i\leq 1$ used before.
\begin{thm}
\label{th2} Assume that $(\xi _i,\mathcal{F}_i)_{i=1,...,n}$ are  martingale differences satisfying (\ref{cdwe}).
Then, for any $ x \geq 0$ and $v > 0$,
\begin{eqnarray}
 \mathbb{P}\left( X_k \geq x\ \mbox{and}\ \langle X\rangle_{k}\leq v^2\ \mbox{for some}\ k \in [1, n]\right)
&\leq& F(x,v) \label{fks} .
\end{eqnarray}
\end{thm}

Bennett \cite{B62} proved (\ref{fks}) in the independent case under the condition that
\[
\mathbb{E}|\xi_{i}|^k \leq \mathbb{E}\xi_{i}^2,\ \ \ \ \ \ \mbox{for any}\ \  k\geq3,
\]
which is in fact equivalent to $|\xi_i|\leq 1$.  Taking into account Remark \ref{remark1}, we see that (\ref{fks}) recovers the inequalities of Freedman and Bennett under the less restrictive condition (\ref{cdwe}).

\section{Proof of Theorems}\label{sec3}
Let $(\xi _i,\mathcal{F}_i)_{i=0,...,n}$ be the supermartingale differences
introduced in the previous section and $X=(X_k,\mathcal{F}_k)_{k=0,...,n}$
be the corresponding supermartingale defined by (\ref{matingal}). For any
nonnegative number $\lambda$, define the \textit{exponential multiplicative
martingale} $Z(\lambda )=(Z_k(\lambda ),\mathcal{F}_k)_{k=0,...,n},$ where
\[
Z_k(\lambda )=\prod_{i=1}^k\frac{e^{\lambda \xi _i}}{\mathbb{E}(e^{\lambda \xi _i}|
\mathcal{F}_{i-1})},\quad \quad  \quad Z_0(\lambda )=1, \ \ \ \ \ \lambda\geq 0.
\]
If $T$ is a stopping time, then $Z_{T\wedge k}(\lambda )$ is also a martingale, where
\[
Z_{T\wedge k}(\lambda )=\prod_{i=1}^{T\wedge k}\frac{e^{\lambda \xi _i}}{\mathbb{E}(e^{\lambda \xi _i}|
\mathcal{F}_{i-1})},\quad \quad  \quad Z_0(\lambda )=1, \ \ \ \ \ \lambda\geq 0.
\]
Thus, for each nonnegative number
$\lambda $ and each $k=1,...,n,$ the random variable $Z_{T\wedge k}(\lambda ) $ is a
probability density on $(\Omega ,\mathcal{F},\mathbb{P})$, i.e.
$$ \int Z_{T\wedge k}(\lambda)  d \mathbb{P} = \mathbb{E}(Z_{T\wedge k}(\lambda))=1.$$ The last observation
allows us to introduce, for any nonnegative number $\lambda ,$ the \textit{%
conjugate probability measure} $\mathbb{P}_\lambda $ on $(\Omega ,\mathcal{F})$
defined by
\begin{equation}
d\mathbb{P}_\lambda =Z_{T\wedge n}(\lambda )d\mathbb{P}.  \label{dp}
\end{equation}
Throughout the paper, we denote by $\mathbb{E}_{\lambda}$ the expectation with
respect to $\mathbb{P}_{\lambda}$.

Consider the predictable process $\Psi (\lambda )=(\Psi _k(\lambda ),%
\mathcal{F}_k)_{k=0,...,n}$, which is called the \textit{cumulant process} and
which is related to the supermartingale $X$ as follows:
\begin{equation}
\Psi _k(\lambda )=\sum_{i=1}^k\log \mathbb{E}(e^{\lambda \xi _i}|\mathcal{F}_{i-1}), \ \ \ 0\leq k \leq n.
\label{C-3}
\end{equation}
We should give a sharp bound for the function $\Psi_{k}(\lambda ).$
To this end, we need the following elementary lemma which, in the special case of centered random variables,
 has been proved by Bennett \cite{B62}.

\begin{lemma}
\label{lemma} If $\xi$ is a random variable such that $\xi \leq 1$, $\mathbb{E}\xi
\leq 0$ and $\mathbb{E}\xi^2=\sigma^2$, then, for any $\lambda \geq 0,$
\begin{equation}\label{fh29}
\mathbb{E}(e^{\lambda \xi } ) \leq \frac{1}{1+\sigma^2} \exp\left\{-\lambda \sigma^2
\right\} + \frac{\sigma^2}{1+\sigma^2 }\exp\{\lambda\} .
\end{equation}
\end{lemma}

\textbf{Proof.}
We argue as in Bennett \cite{B62}. For $\lambda=0$, inequality (\ref{fh29}) is obvious. Fix $\lambda>0$ and consider the function
\begin{eqnarray*}
\phi(\xi) &=& a\xi^2+b\xi +c,\ \ \ \ \xi\leq1,
\end{eqnarray*}
where $a$, $b$ and $c$ are determined by the conditions
\begin{eqnarray*}
\phi(1)=e^{\lambda},\ \ \ \phi(-\sigma^2)=\frac{1}{\lambda}\phi'(-\sigma^2)=\exp\{-\lambda \sigma^2\},\ \ \ \lambda>0.
\end{eqnarray*}
By simple calculations, we have
\begin{eqnarray*}
a&=&\frac{e^{\lambda}-e^{-\lambda\sigma^2}-\lambda(1+\sigma^2)e^{-\lambda%
\sigma^2}}{(1+\sigma^2)^2}, \\
b&=&\frac{ \lambda(1-\sigma^4)e^{-\lambda\sigma^2}+2\sigma^2(e^{%
\lambda}-e^{-\lambda\sigma^2})}{(1+\sigma^2)^2}
\end{eqnarray*}
and
\begin{eqnarray*}
\ \ \ c&=& \frac{\sigma^4e^{\lambda}+(1+2\sigma^2+\lambda\sigma^2+\lambda\sigma^4
) e^{-\lambda\sigma^2} }{(1+\sigma^2)^2}.
\end{eqnarray*}
We now prove that
\begin{eqnarray}
  e^{\lambda \xi} \leq \phi(\xi)  \ \ \ \ \mbox{for any} \ \ \xi \leq1 \ \ \mbox{and}\ \   \lambda>0, \label{keys}
\end{eqnarray}
which will imply the assertion of the lemma.
For any $\xi \in \mathbb{R}$, set $$f(\xi)=\phi(\xi)-e^{\lambda \xi}.$$ Since $f(-\sigma^2)=f(1)=0$, by Rolle's theorem,
there exists some $\xi_1 \in (-\sigma^2, 1)$ such that $f'(\xi_{1})=0$. In the same way, since $f'(-\sigma^2)=0$ and $f'(\xi_{1})=0$,
there exists some $\xi_{2}\in (-\sigma^2, \xi_1)$ such that $f''(\xi_{2})=0$.
Taking into account that the function $f''(\xi)=2a-\lambda^2e^{\lambda \xi}$  is strictly decreasing, we conclude that
$\xi_2$ is the unique zero point of $f''(\xi)$. It follows that $f(\xi)$ is convex on $(-\infty, \xi_{2}]$
and concave on $[\xi_{2}, 1]$, with $\min_{(-\infty, \xi_2]}f(\xi)=f(-\sigma^2)=0$ and $\min_{[\xi_2,1]}f(\xi)=f(1)=0$.
 Therefore $\min_{(-\infty,1]}f(\xi)=0$, which implies  (\ref{keys}).

Since $b\geq0$ and $\mathbb{E}\xi \leq 0$, from (\ref{keys}), it follows that, for any $\lambda>0$,
\begin{eqnarray*}
\mathbb{E}(e^{\lambda \xi}) \leq a\sigma^2 +c = \frac{1}{1+\sigma^2}
\exp\left\{-\lambda \sigma^2 \right\} + \frac{\sigma^2}{1+\sigma^2 }%
\exp\{\lambda\} .
\end{eqnarray*}
This completes the proof of Lemma \ref{lemma}. \hfill\qed

The following technical lemma is from Hoeffding \cite{Ho63} (see Lemma 3 therein and its proof). For reader's convenience,
we shall give a proof following \cite{Ho63}.
\begin{lemma}
\label{lemma2} For any $\lambda \geq0$ and $t\geq 0$, let
\begin{eqnarray}
f(\lambda,t)= \log\left( \frac{1}{1+t}\exp\left\{-\lambda t \right\} + \frac{%
t}{1+t} \exp\{\lambda\} \right). \label{fxt}
\end{eqnarray}
Then $\frac{\partial}{\partial t} f(\lambda, t) >0$ and $\frac{\partial^2}{\partial^2 t} f(\lambda, t) <0$ for any $\lambda >0$ and $t\geq0$.
\end{lemma}

\textbf{Proof.} Denote
\begin{eqnarray*}
g(y) &=& \frac{ e^{-\lambda y}+y-1}{y},\ \ \ \ \ \  y \geq 1.
\end{eqnarray*}
Then $f(\lambda, t)=\lambda+\log g(1+t)$. By straightforward calculation, we have, for any $y\geq 1$,
\[
g'(y)=\frac{e^{-\lambda y}(e^{\lambda y}-1-\lambda y)}{y^2} >0
\]
and
\[
g''(y)=-\frac{2e^{-\lambda y}}{y^3} (e^{\lambda y}-1-\lambda y - \frac{\lambda^2 y^2}{2}) < 0.
\]
Since  $g(y)> 0$ for $y\geq1$,
\[
\frac{\partial}{\partial t} f(\lambda, t)=\frac{g'(y)}{g(y)}\ \ \ \ \ \ \ \ \textrm{and}  \ \ \ \ \ \ \ \    \frac{\partial^2}{\partial^2 t} f(\lambda, t)=\frac{g''(y)g(y)-g'(y)^2}{g(y)^2},
\]
it follows that $\frac{\partial}{\partial t} f(\lambda, t) >0$ and $\frac{\partial^2}{\partial^2 t} f(\lambda, t) <0$ for all $\lambda, t>0$.\hfill\qed

\begin{lemma}
\label{lemma1} Assume that $(\xi _i,\mathcal{F}_i)_{i=1,...,n}$ are supermartingale differences satisfying $\xi_{i} \leq 1$.
Then, for any $\lambda \geq0$ and $k=1,...,n$,
\begin{eqnarray}
\Psi _k(\lambda )  \leq k f\left(\lambda, \frac{\langle X\rangle_{k}}{k} \right). \label{bsf}
\end{eqnarray}
\end{lemma}

\textbf{Proof.} For $\lambda=0$, inequality (\ref{bsf}) is obvious.
By Lemma \ref{lemma}, we have, for any $\lambda > 0$,
\begin{eqnarray}
\mathbb{E}(e^{\lambda \xi _i}|\mathcal{F}_{i-1}) &\leq& \frac{ \exp\left\{-\lambda
\mathbb{E}(\xi_{i}^2|\mathcal{F}_{i-1}) \right\}}{1+\mathbb{E}(\xi_{i}^2|\mathcal{F}_{i-1})} +
\frac{\mathbb{E}(\xi_{i}^2|\mathcal{F}_{i-1})}{1+\mathbb{E}(\xi_{i}^2|\mathcal{F}_{i-1}) }%
\exp\{\lambda\} . \nonumber
\end{eqnarray}
Therefore, using (\ref{fxt}) with $t=\mathbb{E}(\xi_{i}^2|\mathcal{F}_{i-1})$, we get
\begin{eqnarray}
\log\mathbb{E}(e^{\lambda \xi _i}|\mathcal{F}_{i-1}) \leq f(\lambda,\mathbb{E}(\xi_{i}^2|\mathcal{F}
_{i-1})). \label{loge}
\end{eqnarray}
By Lemma \ref{lemma2}, for fixed $\lambda> 0$, the function  $f(\lambda, t)$
has a negative second derivative in $t$.
Hence, $f(\lambda,t)$ is concave in $t\geq 0$, and therefore, by Jensen's inequality,
\begin{eqnarray}
\sum_{i=1}^k f(\lambda,\mathbb{E}(\xi_{i}^2|\mathcal{F}%
_{i-1})) = k \sum_{i=1}^k\frac1k f(\lambda,\mathbb{E}(\xi_{i}^2|\mathcal{F}%
_{i-1}))
 \leq  k f\left(\lambda, \frac{\langle X\rangle_{k}}{k} \right). \label{fjksa}
\end{eqnarray}
Combining  (\ref{loge}) and (\ref{fjksa}), we obtain
\[
\Psi _k(\lambda ) =\sum_{i=1}^k\log\mathbb{E}(e^{\lambda \xi _i}|\mathcal{F}_{i-1})\leq k f\left(\lambda, \frac{\langle X\rangle_{k}}{k} \right).
\]
This completes the proof of Lemma \ref{lemma1}.\hfill\qed

\vspace{0.3cm}

\textbf{Proof of Theorem \ref{th1}.}
 For any $0\leq x \leq n$, define the stopping time
\begin{eqnarray}
 T(x)= \min\{k \in [1,n]:  X_{k}\geq x\ \mbox{and}\ \langle X\rangle_{k}\leq v^2 \},
\end{eqnarray}
with the convention that $\min\{\emptyset\}=0$. Then $$\mathbf{1} \{  X_{k}\geq x\ \mbox{and}\ \langle X\rangle_{k}\leq v^2\ \mbox{for some}\ k \in [1,n]\} =\sum_{k=1}^{n} \mathbf{1} \left\{ T (x)= k\right\}. $$
 Using the change of measure (\ref{dp}), we have, for any  $0\leq x \leq n$, $v>0$ and $\lambda\geq0,$
\begin{eqnarray}
&&\mathbb{P}\left( X_{k}\geq x\ \mbox{and}\ \langle X\rangle_{k}\leq v^2\ \mbox{for some}\ k \in [1,n] \right) \nonumber \\&= & \mathbb{E}_\lambda Z_{T \wedge n}
(\lambda)^{-1}\mathbf{1}_{\left\{  X_{k}\geq x\ \mbox{and}\ \langle X\rangle_{k}\leq v^2\ \mbox{for some}\ k \in[1,n]
\right\}}\nonumber \\
&= & \sum_{k=1}^n \mathbb{E}_\lambda \exp \left\{ -\lambda X_{T \wedge n} +\Psi _{T \wedge n}(\lambda )\right\}
\mathbf{1}_{\left\{ T(x)= k \right\}}  \nonumber \\
&= & \sum_{k=1}^n \mathbb{E}_\lambda \exp \left\{ -\lambda X_k+\Psi _k(\lambda )\right\}\mathbf{1}_{\left\{ T(x)= k \right\}}  \nonumber \\
&\leq & \sum_{k=1}^n \mathbb{E}_\lambda \exp \left\{ -\lambda x +\Psi _k(\lambda )\right\} \mathbf{1}_{\left\{ T(x)= k \right\}}  \label{jhk} .
\end{eqnarray}
Using Lemma \ref{lemma1}, we deduce, for any $0\leq x \leq n$, $v>0$ and $\lambda\geq0,$
\begin{eqnarray}
&&\mathbb{P}( X_k \geq x\ \mbox{and}\ \langle X\rangle_{k}\leq v^2\ \mbox{for some}\ k \in [1,n]  )  \nonumber \\
\ &\leq&\sum_{k=1}^n \mathbb{E}_\lambda \exp \left\{ - \lambda x+k f\left(\lambda, \frac{\langle X \rangle_{k}}{k}\right)  \right\}  \mathbf{1}_{\left\{ T(x)= k \right\}} .
\end{eqnarray}
By Lemma \ref{lemma2}, $f(\lambda, t)$ is increasing in $t\geq0$. Therefore
\begin{eqnarray}
&&\mathbb{P}( X_k \geq x\ \mbox{and}\ \langle X\rangle_{k}\leq v^2\ \mbox{for some}\ k \in [1,n]  )  \nonumber \\
\ &\leq&\sum_{k=1}^n \mathbb{E}_\lambda \exp \left\{ - \lambda x+k f\left(\lambda, \frac{v^2}{k}\right)  \right\}  \mathbf{1}_{\left\{ T(x)= k \right\}} .
\end{eqnarray}
As $f(\lambda,0)=0$  and $f(\lambda,t)$ is concave in $t\geq 0$ (see Lemma \ref{lemma2}), the function  $f(\lambda,t)/t$
is decreasing in $t\geq0$ for any $\lambda\geq 0$.
Hence, we have, for any $0\leq x \leq n$, $v>0$ and $\lambda\geq0,$
\begin{eqnarray}
&&\mathbb{P}( X_k \geq x\ \mbox{and}\ \langle X\rangle_{k}\leq v^2\ \mbox{for some}\ k \in [1,n]  )  \nonumber \\
\ &\leq& \exp \left\{ - \lambda x+n f\left(\lambda, \frac{v^2}{n}\right)  \right\}   \mathbb{E}_\lambda \sum_{k=1}^n \mathbf{1}_{\left\{ T(x)= k \right\}}  \nonumber\\
\ &\leq& \exp \left\{ - \lambda x+n f\left(\lambda, \frac{v^2}{n}\right)  \right\} .   \label{f27}
\end{eqnarray}
Since the function in (\ref{f27}) attains its minimum at
\begin{equation}
\lambda=\lambda(x) = \frac{1}{1+v^2/n}\log\frac{1+x/v^2}{1-x/n } ,
\end{equation}
inserting $\lambda=\lambda(x)$ in (\ref{f27}), we obtain, for any $0\leq x \leq n$ and $v>0$,
\begin{eqnarray}
 \mathbb{P}\left( X_k \geq x\ \mbox{and}\ \langle X\rangle_{k}\leq v^2\ \mbox{for some}\ k \in [1,n] \right)  \leq  H_n(x,v) , \nonumber
\end{eqnarray}
where
\begin{eqnarray}
H_n(x,v)  &=&\inf_{\lambda\geq0} \exp \left\{- \lambda x + n f\left(\lambda, \frac{v^2}n\right) \right\}  ,\label{jksa}
\end{eqnarray}
which gives the bound (\ref{fgl1}).  \hfill\qed

\textbf{Proof of Remark \ref{remark1}.} We will use the function $f(\lambda,t)$ defined by (\ref{fxt}).
Since $\frac{\partial^2}{\partial t^2}f(\lambda,t) \leq 0$ for any $t\geq0$ and $\lambda\geq 0$, it holds
\begin{equation}  \label{flt}
f(\lambda,t)\leq f(\lambda,0)+\frac{\partial }{\partial t }f(\lambda,0)\, t= (e^{\lambda}-1-\lambda)\, t,\ \ \ \ \
t,\lambda\geq0.
\end{equation}
Hence,  using (\ref{jksa}), for
any $x\geq0 $ and $v>0$,
\begin{eqnarray}
H_{n}(x,v)
&\leq&\inf_{\lambda\geq0} \exp \left\{- \lambda x +
(e^{\lambda}-1-\lambda)v^2 \right\}  \nonumber\\
&=& \left(\frac{v^2}{x+v^2}\right)^{x+v^2}e^x ,
\end{eqnarray}
which proves (\ref{freedma1}).
Using the inequality  $$(e^{\lambda}-1-\lambda)v^2\leq \frac{\lambda^2v^2}{2(1-\frac{1}{3}\lambda)},\ \ \ \mbox{for any} \ \lambda, v\geq0,$$ we get, for
any $x\geq0 $ and $v>0$,
\begin{eqnarray}
\left(\frac{v^2}{x+v^2}\right)^{x+v^2}e^x  &\leq&\inf_{3 > \lambda\geq0} \exp \left\{- \lambda x +
\frac{\lambda^2v^2}{2(1-\frac{1}{3}\lambda)} \right\} \nonumber\\
&=& \exp\left\{-\frac{x^2}{v^2\left(1+\sqrt{1+\frac{2\,x}{3\,v^2} } \right)+\frac{1}{3}x }\right\} \nonumber\\
&\leq& \exp\left\{-\frac{x^2}{2(v^2+\frac{1}{3}x )}\right\},   \nonumber
\end{eqnarray}
where the last inequality follows from the fact $\sqrt{1+\frac{2\,x}{3\,v^2} } \leq 1+ \frac{x}{3v^2}$. This proves (\ref{Bennett}) and (\ref{Bernstein}).

Since $f(\lambda,t)/t$ is decreasing in $t\geq0$ for any $\lambda\geq 0$, from (\ref{jksa}), we find that
$H_{n}(x,v)$ is increasing in $n$. Taking into account that $\lim_{n\rightarrow \infty}\left(\frac{n}{n-x}\right)^{n-x}=e^x$, we obtain (\ref{limn}).
This completes the proof of Remark \ref{remark1}.  \hfill\qed

\begin{lemma}
\label{lemma3} Assume that $(\xi _i,\mathcal{F}_i)_{i=1,...,n}$ are martingale differences satisfying
$$\mathbb{E}(\xi_i^2e^{\lambda \xi_{i}}|\mathcal{F}_{i-1}) \leq  e^{\lambda}\mathbb{E}(\xi_i^2|\mathcal{F}_{i-1})$$ for any $\lambda\geq 0$.
Then, for any $\lambda \geq0$ and $k=1,...,n$,
\begin{eqnarray*}
\Psi _k(\lambda ) &\leq& (e^\lambda - 1 -\lambda ) \langle X\rangle_k  .
\end{eqnarray*}
\end{lemma}

\textbf{Proof.} Denote $\psi_i(\lambda)=\log E(e^{\lambda\xi_i}|\mathcal{F}_{i-1})$, $\lambda\geq0$. Since $\psi_i(0)=0$ and $\psi_i'(0)=E(\xi_i|\mathcal{F}_{i-1})=0$, by Leibniz-Newton formula,
it holds
\[
\psi_{i}(\lambda)=\int_{0}^{\lambda}\psi'_{i}(y)dy= \int_{0}^{\lambda} \int_{0}^{y}\psi''_{i}(t)dt dy.
\]
Therefore for any $\lambda \geq0$ and $k=1,...,n$,
\begin{eqnarray}
\Psi _k(\lambda ) = \sum_{i=1}^k\psi_{i}(\lambda)= \sum_{i=1}^k\int_{0}^{\lambda} \int_{0}^{y}\psi''_{i}(t)dt dy . \label{newt}
\end{eqnarray}
 Since, by Jensen's inequality, $\mathbb{E}(e^{t \xi_{i}}|\mathcal{F}_{i-1})\geq 1$, we get, for any $t\geq0$,
\begin{eqnarray}
 \psi''_{i}(t)
&=& \frac{\mathbb{E}(\xi_i^2e^{t \xi_i}|\mathcal{F}_{i-1})}{\mathbb{E}(e^{t \xi_i}|\mathcal{F}_{i-1})} -   \frac{ \mathbb{E}(\xi_i e^{t \xi_i}|\mathcal{F}_{i-1})^2}{ \mathbb{E}(e^{t \xi_i}|\mathcal{F}_{i-1})^2} \nonumber \\
&\leq& \mathbb{E}(\xi_{i}^2e^{t \xi_{i}}|\mathcal{F}_{i-1}) \nonumber\\
&\leq& e^{t}\, \mathbb{E}(\xi_i^2|\mathcal{F}_{i-1}). \label{khjdg}
\end{eqnarray}
Inserting  (\ref{khjdg}) into (\ref{newt}), we obtain
\begin{eqnarray}
\Psi _k(\lambda )
&\leq& \sum_{i=1}^k \int_{0}^{\lambda}\int_{0}^ye^{t}\, \mathbb{E}(\xi_{i}^2|\mathcal{F}_{i-1})dtdy  \nonumber\\
&=& (e^\lambda - 1 -\lambda ) \langle X\rangle_k  \nonumber.
\end{eqnarray}
This completes the proof of Lemma \ref{lemma3}. \hfill\qed

\vspace{0.3cm}

\textbf{Proof of Theorem \ref{th2}.}
By (\ref{jhk}) and Lemma \ref{lemma3}, we obtain, for any $ x\geq 0$, $v>0$ and $\lambda\geq0,$
\begin{eqnarray}
&&\mathbb{P}\left( X_k \geq x\ \mbox{and}\ \langle X\rangle_{k}\leq v^2\ \mbox{for some}\ k \in [1,n]  \right)  \nonumber \\
&\leq & \sum_{k=1}^n \mathbb{E}_\lambda \exp \left\{ -\lambda x +\Psi _k(\lambda )\right\} \mathbf{1}_{\left\{ T(x)= k \right\}}  \nonumber \\
&\leq&\sum_{k=1}^n \mathbb{E}_\lambda \exp \left\{ - \lambda x+(e^\lambda - 1 -\lambda ) \langle X\rangle_k
\right\}  \mathbf{1}_{\left\{ T(x)= k \right\}} \nonumber \\
&\leq&\sum_{k=1}^n \mathbb{E}_\lambda \exp \left\{ - \lambda x+(e^\lambda - 1 -\lambda )v^2
\right\} \mathbf{1}_{\left\{ T(x)= k \right\}} \nonumber \\
&\leq& \exp \left\{- \lambda x + (e^\lambda - 1 -\lambda )v^2 \right\}.\label{dfs}
\end{eqnarray}
Since the function in (\ref{dfs}) attains its minimum at
\begin{equation}
\lambda=\lambda(x) = \log\left(1+\frac{x}{v^2} \right),
\end{equation}
inserting $\lambda=\lambda(x)$ in (\ref{dfs}), we have, for any $ x\geq0$ and $v>0$,
\begin{eqnarray*}
 &&\mathbb{P}\left( X_k \geq x\ \mbox{and}\ \langle X\rangle_{k}\leq v^2\ \mbox{for some}\ k \in [1,n] \right) \\
  &\leq& F(x,v)= \left(\frac{v^2}{x+v^2}\right)^{x+v^2}e^x .
\end{eqnarray*}
This completes the proof of Theorem \ref{th2}.  \hfill\qed

\section{Proof of Corollaries}\label{sec4}

We use Theorem \ref{th1} to prove Corollaries \ref{co2} and \ref{th2.5}.

\vspace{0.2cm}
\textbf{Proof of Corollary \ref{co2}.}
As $-b \leq \xi_{i} \leq 1 $, we have $-\xi_i \leq b$ and $1-\xi_i\geq0$, so that $-\xi_i(1-\xi_i)\leq b(1-\xi_i)$. When $\mathbb{E}(\xi_{i}|\mathcal{F}_{i-1})=0$ or $\mathbb{E}(\xi_{i}|\mathcal{F}_{i-1})\leq0$ and $0< b\leq 1$, it follows that
\begin{eqnarray}
\mathbb{E}(\xi_{i}^2|\mathcal{F}_{i-1}) &=& \mathbb{E}( -\xi_i(1-\xi_i) |\mathcal{F}_{i-1}) + \mathbb{E}(\xi_{i} |\mathcal{F}_{i-1}) \nonumber\\
&\leq & b + (1-b)\mathbb{E}( \xi_i  |\mathcal{F}_{i-1})  \nonumber\\
&\leq&  b. \label{boundf2}
\end{eqnarray}
Therefore $\langle X\rangle_n \leq nb$.
Hence, using Theorem \ref{th1}, we have, for any $0\leq x \leq n$,
\begin{eqnarray}
\mathbb{P}\left(\max_{1\leq k\leq n}X_{k}\geq x \right) &\leq & \sup_{v^2 \leq nb } \mathbb{P}\left(X_k \geq x\ \mbox{and}\ \langle X\rangle_{k}\leq v^2\ \mbox{for some}\ k \in [1, n] \right) \nonumber\\
 &\leq & \sup_{v^2 \leq nb } H_n(x,v) \nonumber \\
 &=&  H_n\left(x,\sqrt{nb}\right),
\end{eqnarray}
which obtains inequality (\ref{Ho11}). Using (\ref{Bernstein}), we get, for any $ x \geq 0$,
\begin{eqnarray}
H_n\left(x,\sqrt{nb}\right) \leq \exp\left\{ -\frac{x^2}{2(nb+\frac13 x)}\right\}. \label{gkns}
\end{eqnarray}
From (\ref{jksa}), we have, for any $ x \geq 0$,
\begin{eqnarray}
 H_n\left(x,\sqrt{nb}\right)
 &=&\inf_{\lambda\geq0} \exp \left\{- \lambda x + n f(\lambda, b)  \right\} . \label{fgky}
\end{eqnarray}
With the notations $z=\lambda(1+b)$ and $p=\frac{b}{1+b}$,  we obtain
 $$f(\lambda, b)=g(z)= - z p + \log(1-p + p e^z).$$
Since $g(0)=g'(0)=0$,
\[
g'(z)=-p+\frac{p}{p+(1-p)e^{-z}}
\]
and
\[
g''(z)=\frac{p(1-p)e^{-z}}{(p+(1-p)e^{-z})^2} \leq \frac{1}{4},
\]
we have
\[
f(\lambda, b)=g(z)\leq \frac{1}{8}z^2= \frac{1}{8}\lambda^2(1+b)^2 .
\]
Returning to (\ref{fgky}),   we deduce, for any $x\geq0$,
\begin{eqnarray}
 H_n\left(x,\sqrt{nb}\right)
 &=&\inf_{\lambda\geq0} \exp \left\{- \lambda x + n f(\lambda, b)  \right\} \nonumber \\
&\leq&\inf_{\lambda\geq0} \exp \left\{- \lambda x + \frac{1}{8} \lambda^2n(1+b)^2  \right\} \nonumber \\
&=&\exp\left\{ - \frac{ 2x^2}{n(1+b)^2}\right\}. \label{fhoes}
\end{eqnarray}
Combining (\ref{gkns}) and (\ref{fhoes}), we obtain (\ref{Ho12}).  \hfill\qed
\vspace{0.3cm}

\textbf{Proof of Corollary \ref{th2.5}.}
For $y >0$ and $k=1,...,n$, set
\begin{eqnarray*}
& X_{k}'=  \sum_{i=1}^{k}  \xi_{i}\textbf{1}_{\{\xi_{i}\leq y \}}, \ \ \ \ \  \  \ \ \ \ \ \ \ &  X_{k}''=\sum_{i=1}^{k}\xi_{i}\textbf{1}_{\{\xi_{i}>y \}}, \\
& X_{k}= X_{k}'+X_{k}'' \ \ \ \ \ \ \ \ \ \ \mbox{and} \ \ \ \ \ \ & V_k^2(y)=\sum_{i=1}^k \mathbb{E}(\xi_{i}^2\mathbf{1}_{\{\xi_{i}\leq y \}} | \mathcal{F}_{i-1}).
\end{eqnarray*}
Since $E(\xi_{i}|\mathcal{F}_{i-1})\leq 0$ implies $E(\xi_{i}\textbf{1}_{\{\xi_{i}\leq y \}}|\mathcal{F}_{i-1})\leq 0$, $X_{k}'$ is a sum of supermartingale differences. Now, for any  $y>0$, $0\leq x \leq ny$ and $v>0$,
\begin{eqnarray*}
&& \mathbb{P}\left(X_k \geq x\ \mbox{and}\ V_k^2(y)\leq v^2 \ \mbox{for some}\ k \in [1,n]\right) \\
&= & \mathbb{P}\left( X_{k}'+X_{k}'' \geq x\ \mbox{and}\ V_k^2(y)\leq v^2 \ \mbox{for some}\ k \in [1,n]\right) \\
 &\leq & \mathbb{P}\left( X_{k}' \geq x\ \mbox{and}\ V_k^2(y)\leq v^2 \ \mbox{for some}\ k \in [1,n] \right) \\
 && +\mathbb{P}\left(  X_{k}'' >0\ \mbox{and}\ V_k^2(y)\leq v^2 \ \mbox{for some}\ k \in [1,n]\right)\\
  &\leq & \mathbb{P}\left( \frac{X_{k}'}{y}  \geq \frac{x}y\ \mbox{and}\ \frac{V_k^2(y)}{y^2} \leq \frac{v^2}{y^2} \ \mbox{for some}\ k \in [1,n]\right) \\
  && + \mathbb{P}\left( \max_{1\leq i \leq n} \xi_{i}>y  \right).
\end{eqnarray*}
Applying Theorem \ref{th1} to $X_k=\frac{X_{k}'}y$, we obtain Corollary \ref{th2.5}.\hfill\qed

\section*{Acknowledgements}
We would like to thank the two referees for their helpful remarks and suggestions.

\end{document}